\documentclass{article}

\usepackage{algorithm, algpseudocode}
\usepackage{amsfonts, amsmath, amsopn, amsthm, epsfig}
  \numberwithin{equation}{section}
\usepackage{graphicx, epstopdf}
\usepackage{hyperref}
\usepackage{subfigure}
\usepackage{appendix}

\usepackage{amssymb}
\usepackage{mathrsfs}

\usepackage{calc}
\usepackage[top=2.5cm,bottom=2.3cm,inner=2.3cm,outer=2.3cm,bindingoffset=0.5cm,footskip=0.55cm]{geometry}
\theoremstyle{remark}

\newenvironment{lemma*}[2][Lemma]{\par\bgroup{\bfseries #1\ #2. }\it\ignorespaces}{\egroup}

\title{SingCubic: Cyclic Incremental Newton-type Gradient Descent with Cubic Regularization for 
Non-Convex Optimization}
\author{Ziqiang Shi\footnotemark[1]
}

\newcommand{\BALD}{\begin{aligned}}
\newcommand{\EALD}{\end{aligned}}
\newcommand{\BALDS}{\begin{aligned*}}
\newcommand{\EALDS}{\end{aligned*}}
\newcommand{\BCAS}{\begin{cases}}
\newcommand{\ECAS}{\end{cases}}
\newcommand{\BEAS}{\begin{eqnarray*}}
\newcommand{\EEAS}{\end{eqnarray*}}
\newcommand{\BEQ}{\begin{equation}}
\newcommand{\EEQ}{\end{equation}}
\newcommand{\BIT}{\begin{itemize}}
\newcommand{\EIT}{\end{itemize}}
\newcommand{\BMAT}{\begin{bmatrix}}
\newcommand{\EMAT}{\end{bmatrix}}
\newcommand{\BNUM}{\begin{enumerate}}
\newcommand{\ENUM}{\end{enumerate}}

\newcommand{\BA}{\begin{array}}
\newcommand{\EA}{\end{array}}

\begin{document}

\maketitle

\renewcommand{\thefootnote}{\fnsymbol{footnote}}

%\footnotetext[1]{Fujitsu Research and Development Center, Beijing, China. E-Mail: shiziqiang@cn.fujitsu.com}

\begin{abstract}
   In this work, we generalized and unified two recent completely different works of~\cite{shi2015large} 
    and~\cite{cartis2012adaptive} respectively into one by proposing the cyclic incremental Newton-type 
    gradient descent with cubic regularization (SingCubic) method for optimizing non-convex functions.
    Through the iterations of SingCubic, a cubic regularized global quadratic 
    approximation using Hessian information is kept and solved.
    Preliminary numerical experiments show the encouraging performance of the SingCubic algorithm when compared to basic incremental or stochastic 
    Newton-type implementations. The results and technique can be served as an initiate for the research on the 
    incremental Newton-type 
    gradient descent methods that employ cubic regularization. The methods and principles proposed in this 
    paper can be used to do logistic regression, autoencoder training, independent components analysis, 
    Ising model/Hopfield network training, multilayer perceptron, deep convolutional network training and so on.
We will open-source parts of our implementations soon.
\end{abstract}

%=============================================
\section{Introduction and Problem Statement}
\label{sec:introduction}

We consider the problem of finding a vector $x^\star\in \mathbb{R}^p$ which minimizes a non-convex 
function\footnotemark[2] $F(x)$, where $F(x)$ is a sum of $n$ sub-functions $f_i (x)$ each is a smooth 
loss function associated with a sample in a training set
\begin{align}
  \min_{x \in \mathbb{R}^p} \,F(x) := \frac{1}{n}\sum_{i=1}^n f_i (x).
  \label{eq:problem}
\end{align}
Problems of this form often arise
in machine learning, such as topic models, dictionary learning, and perhaps most notably, training of deep neural networks.

\footnotetext[2]{Indeed it should be noted that our method proposed in this work also applies to convex functions.}

Since global minimization of non-convex functions is NP-hard~\cite{hillar2013most}, various alternative approaches are applied, for example finding a stationary point or a local minimum. Based on these different approaches previous methods regarding solving~\eqref{eq:problem} can be grouped into two categories also according to what ``order'' of information they require about the objective function.

The most popular of these is the use of iterative optimization methods to reach a stationary point $x^\star$
\begin{align}
  \nabla F(x^\star)=0,
  \label{eq:stationary_point}
\end{align}
for example use stochastic optimization algorithms based on first order information in training deep neural networks.
The standard and popular stochastic gradient descent (SGD)~\cite{robbins1951stochastic} uses iterations of the form
\begin{equation}\label{eqn:sgd}
x_k = x_{k-1} - \alpha_k \nabla f_{i_k}(x_{k-1}),
\end{equation}
where at each iteration an index $i_k$ is sampled uniformly from the set $\{1, ..., n\}$. The randomly
chosen gradient $\nabla f_{i_k}(x_{k-1})$ yields an unbiased estimate of the true gradient $\nabla F(x_{k-1})$ and one can show
under standard assumptions that,  for a suitably chosen decreasing
step-size sequence $\{\alpha_k\}$, the SGD needs $\mathcal{O}(\epsilon^{-4}\text{poly(d)})$  gradient computations of $\nabla f_i(x)$ to find a local $\epsilon$-stationary point that is a point $x$ with $\|\nabla F(x)\|\leq \epsilon$~\cite{ge2015escaping}. The extension work of SGD includes~\cite{ghadimi2016accelerated,allen-zhu2017natasha}. However the stationary point found by these methods is not necessarily the optimal value, it may also be the saddle point or the local maximum.

Besides the first order method, there is another category of methods, called second-order or Newton type methods,
which converge much faster, but need more memory and computation to obtain the curvature information about the objective function. These methods are always employed to escape strict saddle points for non-convex optimization and to find a local minima $x^\star$ satisfying
\begin{align}
  \nabla F(x^\star)=0  \quad \text{and} \quad \nabla^2 F(x^\star)\succeq 0,
  \label{eq:local-optimal-condition}
\end{align}
not necessarily unique. Here strict saddle points are characterized by having negative eigenvalue. In practice the target is always reduced to find a local $\epsilon$-solution
\begin{align}
  \|\nabla F(x^\star)\|\leq \epsilon  \quad \text{and} \quad \lambda_p(\nabla^2 F(x^\star)) \geq -\sqrt{\epsilon},
  \label{eq:epsilon-local-optimal-condition}
\end{align}
where $\lambda_p(H)$ is the smallest eigenvalue of $H$. At a local $\epsilon$-solution, the gradient is guaranteed to be close to zero and the Hessian is guaranteed to be almost positive semidefinite. These second order methods construct a local model of the objective function. Among these algorithms the cubic-regularized Newton's method first considered by~\cite{Griewank1981modification} and more recently by~\cite{nesterov2006cubic} as means for providing the first- and second-order guarantees for the obtained solution. At each iteration, the model used to compute the step from one iterate to the next by solving a sub-problem that approximates the objective function~\eqref{eq:problem} with a cubic-regularized (CR) second-order overestimation at the current iterate $x_k$:
\begin{align}
  \Delta x_k \leftarrow \arg\min_{d \in \mathbb{R}^p} m_k(d)=F(x_k)+d^T\nabla F(x_k)+\frac{1}{2}d^T 
  \nabla^2 F(x_k) d+\frac{M}{6}\|d\|^3.
  \label{eq:cubic-sub-problem}
\end{align}
It has been shown that the CR methods posses the best known iteration complexity to solve~\eqref{eq:problem} within $\mathcal{O}(\epsilon^{-3})$  number of iterations for generating a stationary point and  $\mathcal{O}(\epsilon^{-3/2})$  number of iterations  for generating a local $\epsilon$-solution.

Unfortunately, the CR Newton's method can be
unappealing when $n$ is large or huge, for example in machine learning problems, since its iteration cost scales linearly in $n$.
When the number of components $n$ is very large, then each iteration of~\eqref{eq:cubic-sub-problem} will be very
expensive since it requires computing the gradients and Hessian matrix for all the $n$ component functions $f_i$. More recently~\cite{cartis2012adaptive} presented an inexact cubic-regularized Newton's method (ARC) achieves the same order of theoretical guarantee as the original CR method. ARC relaxed the local model with an inexact Hessian $H_k$ (while the gradient is exact):
\begin{align}
  \Delta x_k \leftarrow \arg\min_{d \in \mathbb{R}^p} m_k(d)=F(x_k)+d^T\nabla F(x_k)+\frac{1}{2}d^T H_k d+\frac{\sigma_k}{3}\|d\|^3,
  \label{eq:arc-sub-problem}
\end{align}
where $H_k$ is sufficiently close to $\nabla^2 F(x_k)$ in the following way
\begin{align}
 \|(H_k-\nabla^2 F(x_k))\Delta x_k\| \leq C\|\Delta x_k\|^2.
  \label{eq:arc-hessian-condition}
\end{align}
In order to further relax the requirement of exact gradient in~\eqref{eq:arc-sub-problem},~\cite{kohler2017sub} propose a pratical sub-sampling scheme (SCR) to implement the inexact cubic-regularized Newton's method:
\begin{align}
  \Delta x_k \leftarrow \arg\min_{d \in \mathbb{R}^p} m_k(d)=F(x_k)+d^Tg_k+\frac{1}{2}d^T H_k d+\frac{\sigma_k}{3}\|d\|^3,
  \label{eq:SCR-sub-problem}
\end{align}
where $g_k=\frac{1}{|S_g(k)|}\sum_{i\in S_g(k)}\nabla f_i(x_k)$ and $H_k=\frac{1}{|S_h(k)|}\sum_{i\in S_h(k)}\nabla^2 f_i(x_k)$, here the two index set $S_g$ and $S_h$ are sampled uniformly from $\{1,...,n\}$ at random (please refer to the complete algorithm in the appendix).

%limited to small-to-medium scale problems that require a high degree of precision.
Based on the related background introduced above, now we can describe our approaches and findings. The primary contribution of this work is the proposal and analysis of a novel algorithm that we call the cyclic incremental Newton-type gradient descent with cubic regularization (SingCubic) method, a cyclic incremental variant of the CR, ARC and SCR method.
The SingCubic method has the low iteration cost as that of SGD methods, but achieves
the convergence rates like the ARC and SCR method stated above. The SingCubic iterations take the form $x_{k+1} \leftarrow x_k + \Delta x_k$, where $\Delta x_k$ is obtained by
\begin{equation}
\label{eqn:singcubic}
    \Delta x_k \leftarrow \arg\min_{d} d^Tg_k+\frac{1}{2}d^T H_k d+\frac{\sigma_k}{3}\|d\|^3,
\end{equation}
where $g_k = \frac{1}{n}\sum_{i=1}^n (g_{k}^i-v_k^i)+H_kx_k$, $H_{k} =  \frac{1}{n}\sum_{i=1}^n H_{k}^i$, and at each iteration, a index $j$ is chosen following a fixed order, and the corresponding $g_{k+1}^j=\nabla f_j(x_{k+1})$, $H_{k+1}^j=\nabla^2 f_j(x_{k+1})$  and $v_{k+1}^j=H_{k+1}^j x_{k+1}$ is selected, then when $i\neq j$ we set
$g_{k+1}^i\leftarrow g_{k}^i$, $H_{k+1}^i\leftarrow H_{k}^i$, and $v_{k+1}^i\leftarrow v_{k}^i$.

That is, like the ARC and SCR methods, the steps incorporates a gradient and a Hessian with respect to each function; but,
like the SGD method, each iteration only computes the gradient and Hessian with respect to a single example (or a single batch of samples) and
the cost of the iterations is independent of $n$.

Besides SingCubic, there are some approaches available to make the CR and ARC methods piratical, and
a full review of this literature would be outside the scope of this work. Several recent work considered randomized variants for stochastic optimization~\cite{kohler2017sub,tripuraneni2017stochastic,ghadimi2017second}. It can be seen that as all these algorithms converge, the required number of gradients and Hessians grows polynomially fast to full data of $n$ samples.
Different from above related methods, the principle behind our SingCubic is similar to that of PROXTONE~\cite{shi2015large}, which kept a global quadratic approximation model of the objective, while in each iteration only use the information of one sample or one batch of samples to update the model. Furthermore the worst-case iteration complexity of SingCubic match those of SCR~\cite{kohler2017sub}, which is the state-of-the-art.
%PROXTONE uses a quadratic model to determine the search direction, which incorporates a user defined positive definite matrix $H_k^j$.

We now outline the rest of the study. Section~\ref{sec:singcubic} presents the main algorithm and gives an equivalent form in order for the ease of analysis. Section~\ref{sec:Analysis} states the assumptions underlying our analysis and gives the main results. We report some experimental results in
Section~\ref{sec:experiments}, and provide concluding remarks in Section~\ref{sec:Conclusions}.

\subsection{Notations and Assumptions}
\label{sec:notations}

In this paper, we assume each $f_i(x)$, for $i=1,\ldots,n$, is differentiable
on the whole space $\mathbb{R}^p$, and their gradients are
Lipschitz continuous,
that is, there exist~$L_i>0$ such that for all $x, y$,
\begin{equation}\label{eqn:smooth-i}
    \|\nabla f_i(x) - \nabla f_i(y)\| \leq L_i \|x-y\|.
\end{equation}

Then from the Lemma 1.2.3 and its proof in Nesterov's book~\cite{Nesterov04book}, for $i=1,\ldots,n$, we have
\begin{equation}\label{eqn:smooth-i-1}
    | f_i(x) -  f_i(y) - \nabla f_i(y)^T(x-y)| \leq \frac{L_i}{2} \|x-y\|^2.
\end{equation}

Similar with the second-order information, it is assumed that the Hessian of each $f_i$ is Lipschitz continuous:
\begin{equation}\label{eqn:second-order-smooth-i}
    \|\nabla^2 f_i(x) - \nabla^2 f_i(y)\| \leq M_i \|x-y\|.
\end{equation}
It is easy to show that
\begin{align}\label{eqn:second-order-upper-bound}
   f_i(x)\leq f_i(y) + \nabla f_i(y)^T(x-y)+ \frac{1}{2}(x-y)^T\nabla^2 f_i(y)(x-y) + \frac{M_i}{6} \|x-y\|^3.
\end{align}

For a symmetric matrix $H$, its spectrum is denoted by ${\lambda_i(H)}_{i=1}^p$. We assume that the eigenvalues
numbered in decreasing order:
\begin{equation}\label{eqn:eigen-in-decreasing-order}
 \lambda_1(H) \geq ... \geq \lambda_p(H).
\end{equation}

%=============================================================
\section{The SingCubic Method}
\label{sec:singcubic}

In this section, we present the novel cyclic incremental Newton-type gradient descent with cubic regularization
(SingCubic) algorithm for solving problems of the form~\eqref{eq:problem}. There are three key steps in the algorithm: (step 2) the cubic-regularized quadratic model~\eqref{eqn:singcubic} is solved to give a search direction; (step 5) update the adaptive weight $\sigma$ of cubic-regularization; (step 6) the component function $f_j(x)$ is sampled cyclically and the cubic regularized quadratic model~\eqref{eqn:singcubic} is updated using this selected function. Once these key steps have been performed, the current point $x_k$ is updated to give a new point of $x_{k+1}$, and the process is repeated.

We summarize the SingCubic method of~\eqref{eqn:singcubic} in Algorithm~\ref{alg:general-singcubic}, while a thorough description of
each of the key steps in the algorithm will follow in the rest of this section.

The SingCubic algorithm may not be understood and reasonable at first glance, however, we can write it in a form as in Algorithm~\ref{alg:general-singcubic-analysis} that is easier to understand. In each iteration, we construct a cubic-regularized quadratic approximation of the original objective. In fact, a quadratic approximation is constructed for each component objective $f_i$ associated with a sample or a batch of samples, however, only one quadratic approximation is updated in each iteration.
To better understand this method, we make the following illustration and observations.

\subsection{The Cubic-Regularized Quadratic Model in Algorithm~\ref{alg:general-singcubic-analysis}}
\label{sec:rqm}

There are two methods in construction the cubic regularized quadratic models.

For fixed $x \in \mathbb{R}^p$ in $k$th iteration, we define a piecewise cubic-regularized quadratic approximation of $F(x)$ as follows:
 \[
F_{k}(x)=\frac{1}{n}\sum_{i=1}^n f_{i,k}(x)
\]
where $f_{i,k}(x)$ is cubic-regularized the quadratic model for $f_i(x)$
\begin{align*}
f_{i,k}(x)=f_i(x_{k_i})+(x-x_{k_i})^T\nabla f_i(x_{k_i})+\frac{1}{2}(x-x_{k_i})^T\nabla^2 f_i(x_{k_i})(x-x_{k_i})+\frac{\sigma_k}{3}\|x-x_{k_i}\|^3,
  %\label{eq:subfunction-surrogate-assum}
\end{align*}
here $k_i$ is a random variable which have the following conditional probability distribution in each iteration:
\begin{align}
\mathbb{P}(k_i=k|j)=\frac{1}{n} \quad \text{and} \quad \mathbb{P}(k_i=k_{i-1}|j)=1-\frac{1}{n},
  \label{eq:random-parameter-dist}
\end{align}
Then at each iteration the search direction is found by solving the subproblem~\eqref{eq:search-direction}.

\begin{algorithm}[H]
\caption{SingCubic: A generic stochastic incremental Newton-type gradient descent with cubic regularization}
\label{alg:general-singcubic}

\textbf{Input and initialization}: Start point $x_0 \in$ dom $F$; let $k=0$, for $i\in\{1,2,..,n\}$, $g_k^i=\nabla f_i(x_k)$, $H_k^i=\nabla^2 f_i(x_k)$, $f_k^i=f_i(x_k)$, $u_k^i=g_k^{iT}x_k$, $v_k^i=H_k^i x_k$, $w_k^i=x_k^TH_k^i x_k$; let $H_{k} =  \frac{1}{n}\sum_{i=1}^n H_{k}^i$, $g_k = \frac{1}{n}\sum_{i=1}^n (g_{k}^i-v_k^i)+H_kx_k$, $c_k = g_k^Tx_k-\frac{1}{2}x_k^TH_kx_k+\frac{1}{n}\sum_{i=1}^n \left[f_k^i- u_k^i+\frac{1}{2}w_k^i\right]$; and $\eta_1=0.1$, $\eta_2=0.9$, $\gamma_1=2.0$, $\gamma_2=2.0$, $\sigma=1.0$, $\epsilon_g=1e-6$, $\lambda_k=0$.

1: \textbf{Repeat}

2:  \quad   Solve the subproblem for a search direction:
\begin{align*}
\Delta x_k \leftarrow \arg\min_{d} m_k(d)=c_k+d^Tg_k+\frac{1}{2}d^T H_k d+\frac{\sigma_k}{3}\|d\|^3.
  \label{eq:subproblem_prime}
\end{align*}

3: \quad Compute $F(x_k+\Delta x_k)$ and
\begin{align}
\rho_k=\frac{F(x_k)-F(x_k+\Delta x_k)}{F(x_k)-m_k(\Delta x_k)}.
\end{align}

4: \quad If $\rho_k >= \eta_1$, update $x_{k+1}=x_k+\Delta x_k$; otherwise $x_{k+1}=x_k$.

5: \quad If $\rho_k >= \eta_2$, $\sigma_{k+1}=\max(\sigma_k/\gamma_2,1e-16)$; else if $\rho_k < \eta_1$, $\sigma_{k+1}=\gamma_1*\sigma_k$.

6: \quad Sample a  index $j$ fowlling a fixed order and corresponding $g_{k+1}^j=\nabla f_j(x_{k+1})$, $H_{k+1}^j=\nabla^2 f_j(x_{k+1})$, $u_{k+1}^j=g_{k+1}^Tx_{k+1}$, $v_{k+1}^j=H_{k+1}^j x_{k+1}$, $w_{k+1}^j=x_{k+1}^TH_{k+1}^j x_{k+1}$ is selected, then when $i\neq j$ we update
$g_{k+1}^i\leftarrow g_{k}^i$, $H_{k+1}^i\leftarrow H_{k}^i$,  $u_{k+1}^i\leftarrow u_{k}^i$,  $v_{k+1}^i\leftarrow v_{k}^i$, and $w_{k+1}^i\leftarrow w_{k}^i$; further $H_k$, $g_k$, and $c_k$ is updated as

\begin{eqnarray*}
   H_{k+1}&=&\frac{1}{n}\sum_{i=1}^n H_{k+1}^j, \\
   g_{k+1}&=&\frac{1}{n}\sum_{i=1}^n \left[-v_{k+1}^i + g_{k+1}^i \right]+H_{k+1} x_{k+1},\\
   c_{k+1}&=&g_{k+1}^Tx_{k+1}-\frac{1}{2}x_{k+1}^TH_{k+1}x_{k+1}+\frac{1}{n}\sum_{i=1}^n \left[f_{k+1}^i- u_{k+1}^i+\frac{1}{2}w_{k+1}^i\right].
\end{eqnarray*}

7: \quad $k \leftarrow k+1$.

8: \textbf{Until} stopping conditions are satisfied.

\textbf{Output}: $x_k$.
\end{algorithm}

We can express the derivatives of the cubic-regularized the quadratic model $f_{i,k}(x)$ as 
\begin{align}
\nabla_x f_{i,k}(x) = \nabla f_i(x_{k_i})+\nabla^2 f_i(x_{k_i})(x-x_{k_i})+\sigma_k\|x-x_{k_i}\|(x-x_{k_i}),
    \label{eq:grad-subfunction-surrogate-assum}
  \end{align}
and
\begin{align}
\nabla_{xx}f_{i,k}(x) =\nabla^2 f_i(x_{k_i})+\sigma_k\|x-x_{k_i}\|I+\sigma_k\|x-x_{k_i}\|
  (\frac{x-x_{k_i}}{\|x-x_{k_i}\|})(\frac{x-x_{k_i}}{\|x-x_{k_i}\|})^T,
    \label{eq:hess-subfunction-surrogate-assum}
  \end{align}

One of the crucial ideas of this algorithm is that the component function to be used for updating the search direction
at each iteration is chosen randomly. This allows the function to be selected very
quickly. After the component function $f_j(x)$ selected and updated by~\eqref{eq:subfunction-surrogate-update},
while leaving all other $f_{j,k+1}(x)$ unchanged.

\begin{algorithm}[H]
\caption{SingCubic in an equivalent form}
\label{alg:general-singcubic-analysis}

\textbf{Input and initialization}: Start point $x_0 \in$ dom $f$; $k=0$, for $i\in\{1,2,..,n\}$, let $f_{i,k}(x)=f_i(x_k)+(x-x_k)^Tg_{k}^i+\frac{1}{2}(x-x_k)^TH_{k}^i(x-x_0)$, where $g_k^i=\nabla f_i(x_k)$, $H_k^i=\nabla^2 f_i(x_k)$; and $F_k(x)=\frac{1}{n}\sum_{i=1}^n f_{i,k}(x)$; and $\eta_1=0.1$, $\eta_2=0.9$, $\gamma_1=2.0$, $\gamma_2=2.0$, $\sigma=1.0$, $\epsilon_g=1e-6$, $\lambda_k=0$.

1: \textbf{Repeat}

2: \quad Solve the subproblem for new approximation of the solution:
\begin{align}
\Delta x_k \leftarrow \arg\min_{x} \bigl[ m_k(d)=F_k(x_k+d) + \frac{\sigma_k}{3}\|d\|^3 \bigr] .
  \label{eq:search-direction}
\end{align}

3: \quad Compute $F(x_k+\Delta x_k)$ and
\begin{align}
\rho_k=\frac{F(x_k)-F(x_k+\Delta x_k)}{F(x_k)-m_k(\Delta x_k)}.
\end{align}

4: \quad If $\rho_k >= \eta_1$, update $x_{k+1}=x_k+\Delta x_k$; otherwise $x_{k+1}=x_k$.

5: \quad If $\rho_k >= \eta_2$, $\sigma_{k+1}=\max(\sigma_k/\gamma_2,1e-16)$; else if $\rho_k < \eta_1$, $\sigma_{k+1}=\gamma_1*\sigma_k$.

6: \quad Sample $j$ from $\{1,2,..,n\}$ in a fixed order, and update the quadratic models:
\begin{align}
f_{j,k+1}(x)=f_j(x_{k+1})+(x-x_{k+1})^T\nabla f_j(x_{k+1})+\frac{1}{2}(x-x_{k+1})^TH_{i,k+1}(x-x_{k+1}),
  \label{eq:subfunction-surrogate-update}
\end{align}
while leaving all other $f_{i,k+1}(x)$ unchanged: $f_{i,k+1}(x)\leftarrow f_{i, k}(x)$ ($i\neq j$); and $F_{k+1}(x)=\frac{1}{n}\sum_{i=1}^n f_{i,k+1}(x)$.

7: \textbf{Until} stopping conditions are satisfied.

\textbf{Output}: $x_k$.
\end{algorithm}

\subsection{The Subproblem}

The subproblem~(\ref{eq:search-direction}) or~(\ref{eq:subproblem_prime}) is a unconstrained cubic regularized quadratic optimization problem, and
several methods have been proposed to solve it more efficiently.~\cite{cartis2012adaptive} proposed to
approximately solve it in Krylov space.~\cite{agarwal2017finding} proposed an alternative fast way to
solve it.~\cite{carmon2016gradient} proposed
a method based on gradient descent.
In order for completeness, an exact solver based on trust region methods~\cite{conn2000trust} is included in this paper, and 
it is summarized in Algorithm~\ref{alg:prox_lasso}.

\begin{align}
x^{k+1} \leftarrow \arg\min_{x} \bigl[ G^k(x) + \lambda_2\|x\|_1 \bigr]\nonumber\\ = \arg\min_{x} F^k(x).
  \label{eq:search-direction-l1lr}
\end{align}

That means for each gradient, we need to use several iterations of computing approximated Hessian to forming a lasso problem, which also needs several iterations to solve. Thus typically SingCubic needs much more time for each iteration than that of the first-order method. 

%That means although SingCubic is much fast than other methods in the number of gradients or epochs, but maybe slower in time.

\subsection{The Weight $\sigma$ of Cubic-Regularization}
\label{sec:sigma}

 As the Algorithm~\ref{alg:general-singcubic-analysis} shows, $\sigma$ changes following the similar principle as~\cite{cartis2012adaptive}, that is, if the cubic-regularized the quadratic approximation model is very close to the original function,  $\sigma$ will not change or decrease; if the approximation model is far away from the original function, it will increase.

\subsection{The Update of $g$, $H$, $c$}
\label{sec:updateghc}

The derivation of SingCubic:

Let
\begin{align}
c_k+(x-x_k)^Tg_k+\frac{1}{2}(x-x_k)^T H_k(x-x_k)=\frac{1}{n}\sum_{i=1}^nf_{i,k}(x),
\label{eq:equavelence}
\end{align}
where
\begin{align}
f_{i,k}(x) =f_i(x_{k_i})+(x-x_{k_i})^T\nabla f_i(x_{k_i})+\frac{1}{2}(x-x_{k_i})^T\nabla^2 f_i(x_{k_i})(x-x_{k_i}).
\end{align}
Find the second derivative with respect to $x$ on both sides of~\eqref{eq:equavelence}, we have
\begin{align}
H_k=\frac{1}{n}\sum_{i=1}^n\nabla^2 f_i(x_{k_i}).
\end{align}
Take the first derivative with respect to $x$ on both sides of~\eqref{eq:equavelence}, we have
\begin{align}
g_k+H_k(x-x_k)=\frac{1}{n}\sum_{i=1}^n \left[\nabla^2 f_i(x_{k_i})(x-x_{k_i}) + \nabla f_i(x_{k_i}) \right],
\end{align}
and set $x=x_k$, that is
\begin{align}
g_k=\frac{1}{n}\sum_{i=1}^n \left[\nabla^2 f_i(x_{k_i})(x_k-x_{k_i}) + \nabla f_i(x_{k_i}) \right].
\end{align}
Do some simplification we have
\begin{align}
g_k=\frac{1}{n}\sum_{i=1}^n \left[-\nabla^2 f_i(x_{k_i}) x_{k_i}+ \nabla f_i(x_{k_i}) \right]+H_k x_k.\nonumber  
\end{align}
Set $x=x_k$ directly at both sides of the expression~\eqref{eq:equavelence}, we have
\begin{align}
c_k=\frac{1}{n}\sum_{i=1}^nf_{i,k}(x_k)=\frac{1}{n}\sum_{i=1}^n [f_i(x_{k_i})+(x_k-x_{k_i})^T\nabla f_i(x_{k_i})+\frac{1}{2}(x_k-x_{k_i})^T\nabla^2 f_i(x_{k_i})(x_k-x_{k_i})],
\end{align}
that is
\begin{align}
c_k=g_k^Tx_k-\frac{1}{2}x_k^TH_kx_k+\frac{1}{n}\sum_{i=1}^n [f_i(x_{k_i})-\nabla f_i(x_{k_i})^Tx_{k_i}+\frac{1}{2}x_{k_i}^T\nabla^2 f_i(x_{k_i})x_{k_i}].
\end{align}

  \begin{algorithm}[H]
  \caption{Solving subproblem~(\ref{eq:search-direction})}
  \label{alg:prox_lasso}
  
  \textbf{Input}:  $k=0$, starting point $g\in \mathbb{R}^p$, $H$,
 $\lambda_k$,  $\epsilon_{tol}=0.1$, $\epsilon=2\sqrt{\text{2.22e-16}}$.

  1: Compute the bounds of $\lambda$.

  \quad 1.1: By using Gershgorin Circle Theorem,  the lower and upper bounds of the eigenvalue can be 
  computed as
  \begin{align*}
    G_l \leftarrow \min_{i={1,...,p}} (H(i,i)-\sum_{j=1,j\neq i}^p|H(i,j)|),
    \end{align*}
and
    \begin{align*}
      G_u \leftarrow \max_{i={1,...,p}} (H(i,i)+\sum_{j=1,j\neq i}^p|H(i,j)|),
      \end{align*}
respectively.

\quad 1.2: Find the larger root $\lambda_1$  of $x^2+G_lx-\|g\|\sigma=0$; 
find the larger root $\lambda_{2}$ of $x^2+G_ux-\|g\|\sigma=0$.

\quad 1.3: Set $\lambda_{lower}=\max(0,-\min_{i=1,...,p}(H(i,i)),\lambda_{2})$ and 
$\lambda_{upper}=\max(0,\lambda_{1})$.

2: Reinitialize at previous lambda in case of unsuccessful iterations:
  If $\lambda_{lower}\leq \lambda_k \leq \lambda_{upper}$, $\lambda= \lambda_k$; else $\lambda$ 
  is uniformly random sampling from $[\lambda_{lower}, \lambda_{upper}]$.

3: \textbf{Repeat}

\quad 3.1: Let $H(\lambda)= H+\lambda I$. If $\lambda_{lower}==\lambda_{upper}==0$ or 
 $g==0$, then $\lambda\in \mathcal{N}$ and goto 3.8.

\quad 3.2: Do the Cholesky factorization of $H(\lambda)=LL^T$, If factorization is 
failed, then $\lambda\in \mathcal{N}$ and goto 3.8.

\quad 3.3:  Solve $LL^Td=-g$ for $d$.

\quad 3.4:  Let $\phi(\lambda)=\frac{1}{\|d\|}-\frac{\sigma}{\lambda}$, 
if $|\phi(\lambda)|\leq\epsilon_{tol}$, 
then break and \textbf{return} $d$.

\quad 3.5: Let $w=L^{-1}d$.

\quad 3.6: If $\phi(\lambda)<0$, then $\lambda\in \mathcal{L}$, let $\lambda_{lower}\leftarrow\lambda$, 
compute  $c_{hi}$ as the max root of 
$\frac{\|w\|^2}{\|d\|^3}x^2+(\frac{1}{\|d\|}+\frac{\|w\|^2}{\|d\|^3}\lambda)x+
\frac{1}{\|d\|}\lambda-\sigma=0$, 
set $\lambda^+\leftarrow \lambda+c_{hi}$, and $\lambda\leftarrow \lambda^+$; 

\quad  3.7: Else if $\phi(\lambda)>0$, then $\lambda\in \mathcal{G}$,
let $\lambda_{upper}\leftarrow\lambda$, compute  $c_{hi}$ as the max root of 
$\frac{\|w\|^2}{\|d\|^3}x^2+(\frac{1}{\|d\|}+\frac{\|w\|^2}{\|d\|^3}\lambda)x+
\frac{1}{\|d\|}\lambda-\sigma=0$, and set $\lambda^+\leftarrow\lambda+c_{hi}$.

\quad  \quad  3.7.1: If $\lambda^+>0$, do the Cholesky factorization of $H+\lambda^+I=LL^T$,
and set $\lambda\leftarrow \lambda^+$, If factorization is 
failed, then $\lambda^+\in \mathcal{N}$.

\quad \quad  3.7.2: If $\lambda^+\leq 0$ or $\lambda^+\in \mathcal{N}$, $\lambda_{lower}\leftarrow \max(\lambda_{lower}, 
\lambda^+)$, $\lambda=\max(\sqrt{\lambda_{lower}*\lambda_{upper}},\lambda_{lower}+
0.01*(\lambda_{upper}-\lambda_{lower}))$, then if $\lambda_{upper}==\lambda_{lower}$, 
$\lambda\leftarrow \lambda_{lower}$, do the eigen decomposition of $H$, 
get the eigenvalues $\Lambda$ 
and eigenvectors $U$, 
   let $u_p$ be the eigenvector corresponding to the smallest eigenvalue; 
   and compute $\alpha$ as the min root of 
   $x^2 + 2u_p^Tdx + d*d-\frac{\lambda^2}{\sigma^2}=0$,
  set $d\leftarrow d+\alpha * u_p$; then break and \textbf{return} $d$.

  \quad 3.8: If  $\lambda\in \mathcal{N}$, do the following
  
  \quad \quad 3.8.1: $\lambda_{lower}\leftarrow \max(\lambda_{lower}, 
  \lambda)$, $\lambda\leftarrow \max(\sqrt{\lambda_{lower}*\lambda_{upper}},\lambda_{lower}+
  0.01*(\lambda_{upper}-\lambda_{lower}))$. 

  \quad \quad 3.8.2: If $\lambda_{upper}==\lambda_{lower}$, do

  \quad \quad \quad 3.8.2.1: $\lambda\leftarrow \lambda_{lower}$, do the eigen decomposition of $H$, 
  get the eigenvalues $\Lambda$ 
  and eigenvectors $U$.

  \quad \quad \quad 3.8.2.2: If all the the eigenvalues $\Lambda\geq 0$, then break and \textbf{return} $d$.

  \quad \quad \quad 3.8.2.3:  Let $u_p$ be the eigenvector corresponding to the smallest eigenvalue; 
     and compute $\alpha$ as the min root of 
     $x^2 + 2u_p^Tdx + d*d-\frac{\lambda^2}{\sigma^2}=0$,
    set $d\leftarrow d+\alpha * u_p$; then break and \textbf{return} $d$.

\quad \textbf{Until}  max number 
of iterations is achieved.

  \textbf{Output}: $d$.
  \end{algorithm}

\subsection{The Inexact SingCubic Method}
\label{sec:isingcubic}

Each iteration of the SingCubic method requires the calculation of the special cubical proximity operator
\begin{align}
\text{cuprox}_{H,M}(y) = \arg\min_{x} \bigl[ \frac{1}{2}\|x-y\|_H^2 + \frac{M}{6}\|x-y\|^3 \bigr] .
  \label{eq:cuprox}
\end{align}
This cuprox operator has no analytic solution, and it is very expensive to compute the solution exactly. Despite the difficulty in computing the exact cuprox operator, efficient methods have been developed to compute approximate the cuprox operator. In this work, we show in several contexts that, provided the error in the cuprox operator calculation is controlled in an appropriate way, inexact SingCubic achieve the same convergence rates as the exact SingCubic algorithm.

We use $\varepsilon_k$ to denote the error in the calculation of cuprox operator achieved by $x_k$, meaning that
\begin{align}
\frac{1}{2}\|x_k-y\|_H^2 + \frac{M}{6}\|x_k-y\|^3 \leq \varepsilon_k +\min_{x} \bigl[ \frac{1}{2}\|x-y\|_H^2 + \frac{M}{6}\|x-y\|^3 \bigr] .
  \label{eq:ecuprox}
\end{align}

The inexact version of SingCubic algorithm is obtained by simply replacing (\ref{eq:search-direction}) by finding a $x_{k+1}$ satisfying
\begin{align}
F_k(x_{k+1}) + \frac{M}{6}\|x_{k+1}-x_k\|^3 \leq \varepsilon_k + \min_{x} \bigl[ F_k(x) + \frac{M}{6}\|x-x_k\|^3 \bigr] .
  \label{eq:inexact-search-direction}
\end{align}

\subsection{Convergence Analysis}
\label{sec:Analysis}

The techniques provided in~\cite{kohler2017sub,shi2015large} can be smoothly adapted to prove the convergence of SingCubic, and theoretically has the same convergence speed as SCR in~\cite{kohler2017sub}. But in experiments, it can be seen that SingCubic is faster and it is more efficient in using gradients.

\begin{figure*}[th]%[htp]
%\vspace{-0.4in}
\centering
%\begin{center}
\begin{tabular}{ccc}
\hspace{-5mm}
\begin{tabular}{c}
(a)\includegraphics[width=0.45\linewidth]{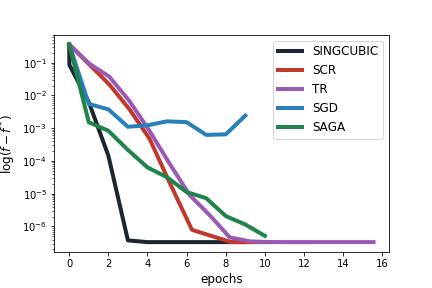}
\end{tabular}
\hspace{-5mm}
 &
\begin{tabular}{c}
(b)\includegraphics[width=0.45\linewidth]{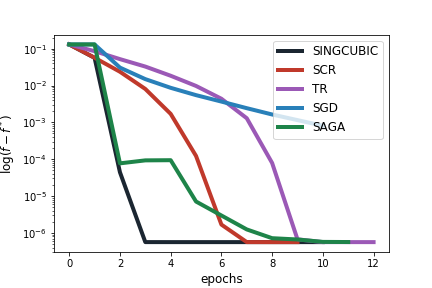}
\end{tabular}
\hspace{-5mm}
\end{tabular}
%\end{center}
%\vspace{-2.5mm}
\caption{
A comparison of SingCubic to competing optimization techniques in solving convex logistic regression for two datasets, (a) is a9a; (b) is covertype.
}
\label{fig_results}
%\vspace{-2.5mm}
\end{figure*}

\section{Experimental Results}
\label{sec:experiments}

Here present the 
results of some numerical experiments to illustrate the properties of the SingCubic method. 

Five algorithms are evaluated and compared:
\begin{itemize}
\item
SingCubic: Penalty increase multiplier is 2.0, 
penalty decrease multiplier is 2.0, the initial penalty parameter is 0.01, initial tr radius is 1.0, 
the successful threshold is 0.1, the very successful threshold is 0.9, the batch size is set as 0.001*$n$.
\item
SCR: Penalty increase multiplier is 2.0, 
penalty decrease multiplier is 2.0, the initial penalty parameter is 0.01, initial tr radius is 1.0, 
the successful threshold is 0.1, the very successful threshold is 0.9. 
\item
TR: Penalty increase multiplier is 2.0, 
penalty decrease multiplier is 2.0, the initial penalty parameter is 0.01, initial tr radius is 1.0, 
the successful threshold is 0.1, the very successful threshold is 0.9. 
\item
SGD: Learning rate 0.1, batch size is set as 0.001*$n$.
\item
SAGA: Learning rate 0.01, initial weights are set to 0.
\end{itemize}

\begin{figure*}[th]%[htp]
%\vspace{-0.4in}
\centering
%\begin{center}
\begin{tabular}{ccc}
\hspace{-5mm}
\begin{tabular}{c}
(a)\includegraphics[width=0.45\linewidth]{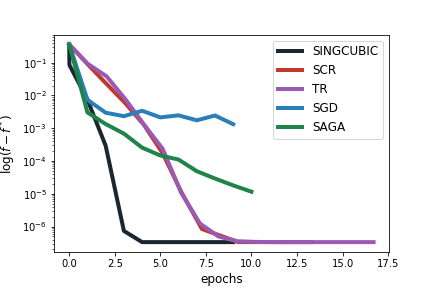}
\end{tabular}
\hspace{-5mm}
 &
\begin{tabular}{c}
(b)\includegraphics[width=0.45\linewidth]{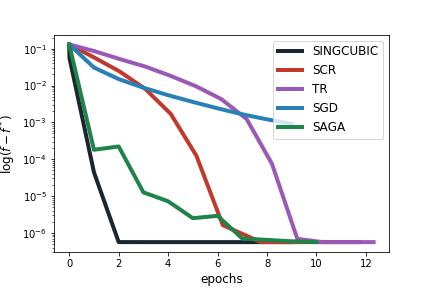}
\end{tabular}
\hspace{-5mm}
\end{tabular}
%\end{center}
%\vspace{-2.5mm}
\caption{
A comparison of SingCubic to competing optimization techniques in solving nonconvex logistic regression for two datasets, (a) is a9a; (b) is covertype.
}
\label{fig_results_nonconvex}
%\vspace{-2.5mm}
\end{figure*}

We focus on 
the convex/nonconvex logistic regression problem for binary classification: given a set of training 
examples
$(x_1,y_1),\ldots, (x_n, y_n)$ where $x_i\in\mathbb{R}^p$ and $y_i\in\{0, 1\}$, we find the optimal 
predictor $w\in\mathbb{R}^p$ by solving the convex $L2$-regularized convex problem
\begin{align}
\min_{w\in\mathbb{R}^p} \quad \,F(w) :=
    -\frac{1}{n} \sum_{i=1}^n \bigl[ y_i\log\frac{1}{1+\exp(-x_i^T w)}+(1-y_i)\log(1-\frac{1}{1+\exp(-x_i^T w)})\bigr]
    + \alpha \frac{1}{2}\|w\|_2^2,
  \label{eq:logistic-loss}
\end{align}
or the nonconvex  logistic regression
\begin{align}
  \min_{w\in\mathbb{R}^p} \quad \,F(w) := 
      -\frac{1}{n} \sum_{i=1}^n \bigl[ y_i\log\frac{1}{1+\exp(-x_i^T w)}+(1-y_i)\log(1-\frac{1}{1+\exp(-x_i^T w)})\bigr]
    + \alpha \sum_{j=1}^p \frac{\beta w_j^2}{1+\beta w_j^2},
    \label{eq:nonconvex-logistic-loss}
  \end{align}
where $\alpha$ and $\alpha$ are the regularization parameters.

The $X$ and $Y$ are from the popular a9a and covertype datasets. For a9a, $n=32561$, $p=123$, and $Y\in \{0,1\}$; for covertype, $n=581012$, $d=54$, and $Y\in \{1,2\}$.

The results of the different methods are plotted for the effective epochs for a9a and covertype respectively through the data in Figure~\ref{fig_results} and Figure~\ref{fig_results_nonconvex} for convex and nonconvex problems respectively. The iterations of SingCubic seem to achieve the best of all.

%\subsection{multinominal (softmax) regression on  mnist}
%\label{sec:multreg}

%\subsection{Rosenbrock}
%\label{sec:rosenbrock}

\section{Conclusions}
\label{sec:Conclusions}

This paper introduces a novel cyclic incremental Newton-type 
    gradient descent with cubic regularization method called SingCubic for minimizing non-convex finite sums. 
In the algorithm description and empirical study, we make clear the implementation details of SingCubic and do the numerical evaluations to both convex and nonconvex problems. 
We show that SingCubic can make full use of gradients and Hessians, converges much faster than state-of-the-art second-order methods in the number of epochs.

% ---- Bibliography ----
%\bibliographystyle{abbrv}  %this one
%\bibliography{IEEEabrv,CommunityBIB-Jerry.bib}
\bibliographystyle{splncs03}
\bibliography{singcubic}

\begin{thebibliography}{10}
\providecommand{\url}[1]{\texttt{#1}}
\providecommand{\urlprefix}{URL }

\bibitem{agarwal2017finding}
{Agarwal}, N., {Allen-Zhu}, Z., {Bullins}, B., {Hazan}, E., {Ma}, T.: Finding
  approximate local minima faster than gradient descent. In: Proceedings of the
  49th Annual ACM SIGACT Symposium on Theory of Computing. pp. 1195--1199
  (2017)

\bibitem{allen-zhu2017natasha}
{Allen-Zhu}, Z.: Natasha 2: Faster non-convex optimization than sgd. arXiv
  preprint arXiv:1708.08694  (2017)

\bibitem{carmon2016gradient}
{Carmon}, Y., {Duchi}, J.C.: Gradient descent efficiently finds the
  cubic-regularized non-convex newton step. arXiv preprint arXiv:1612.00547
  (2016)

\bibitem{cartis2012adaptive}
{Cartis}, C., {Gould}, N.I.M., {Toint}, P.L.: An adaptive cubic regularization
  algorithm for nonconvex optimization with convex constraints and its
  function-evaluation complexity. Ima Journal of Numerical Analysis  32(4),
  1662--1695 (2012)

\bibitem{conn2000trust}
Conn, A.R., Gould, N.I., Toint, P.L.: Trust region methods, vol.~1. Siam (2000)

\bibitem{ge2015escaping}
{Ge}, R., {Huang}, F., {Jin}, C., {Yuan}, Y.: Escaping from saddle points ---
  online stochastic gradient for tensor decomposition. conference on learning
  theory pp. 797--842 (2015)

\bibitem{ghadimi2016accelerated}
{Ghadimi}, S., {Lan}, G.: Accelerated gradient methods for nonconvex nonlinear
  and stochastic programming. Mathematical Programming  156,  59--99 (2016)

\bibitem{ghadimi2017second}
{Ghadimi}, S., {Liu}, H., {Zhang}, T.: Second-order methods with cubic
  regularization under inexact information  (2017)

\bibitem{Griewank1981modification}
{Griewank}, A.: The modification of newton's method for unconstrained
  optimization by bounding cubic terms. Technical report NA/12  (1981)

\bibitem{hillar2013most}
{Hillar}, C.J., {Lim}, L.H.: Most tensor problems are np-hard. Journal of the
  ACM  60(6), ~45 (2013)

\bibitem{kohler2017sub}
{Kohler}, J.M., {Lucchi}, A.: Sub-sampled cubic regularization for non-convex
  optimization. international conference on machine learning  70,  1895--1904
  (2017)

\bibitem{Nesterov04book}
Nesterov, Y.: Introductory Lectures on Convex Optimization: A Basic Course.
  Kluwer, Boston (2004)

\bibitem{nesterov2006cubic}
{Nesterov}, Y., {Polyak}, B.T.: Cubic regularization of newton method and its
  global performance. Mathematical Programming  108(1),  177--205 (2006)

\bibitem{robbins1951stochastic}
{Robbins}, H., {Monro}, S.: A stochastic approximation method. Annals of
  Mathematical Statistics  22(3),  400--407 (1951)

\bibitem{shi2015large}
{Shi}, Z., {Liu}, R.: Large scale optimization with proximal stochastic
  newton-type gradient descent. In: ECMLPKDD'15 Proceedings of the 2015th
  European Conference on Machine Learning and Knowledge Discovery in Databases
  - Volume Part I. pp. 691--704 (2015)

\bibitem{tripuraneni2017stochastic}
{Tripuraneni}, N., {Stern}, M., {Jin}, C., {Regier}, J., {Jordan}, M.I.:
  Stochastic cubic regularization for fast nonconvex optimization. arXiv
  preprint arXiv:1711.02838  (2017)

\end{thebibliography}

\end{document}